\documentclass[12pt]{article}
\usepackage{amssymb}
\usepackage{amsmath}
\usepackage{amsthm}
\usepackage{amsbsy}
\usepackage{graphicx}
\usepackage[]{epstopdf}
\newcommand{\be}{\begin{equation}}
\newcommand{\ee}{\end{equation}}
\newcommand{\bea}{\begin{eqnarray}}
\newcommand{\eea}{\end{eqnarray}}
\newcommand{\ba}{\begin{array}}
\newcommand{\ea}{\end{array}}

\newcommand{\bc}{\begin{center}}
\newcommand{\ec}{\end{center}}
\newcommand{\ben}{\begin{enumerate}}
\newcommand{\een}{\end{enumerate}}
\newcommand{\bfi}{\begin{figure}}
\newcommand{\efi}{\end{figure}}

\newcommand{\bq}{\begin{quote}}
\newcommand{\eq}{\end{quote}}
\newcommand{\bqu}{\begin{quotation}}
\newcommand{\equ}{\end{quotation}}
\newenvironment{emphit}{\begin{itemize}}{\end{itemize}}
\newcommand{\bemp}{\begin{emphit}}
\newcommand{\eemp}{\end{emphit}}

\newcommand{\bt}{\begin{tabular}}
\newcommand{\et}{\end{tabular}}

\newtheorem{myth}{Theorem}[section]
\newtheorem{mylem}{Lemma}[section]

\newtheorem{mypro}{Proposition}[section]
\newtheorem{mydef}{Definition}[section]
\newtheorem{myrem}{Remark}[section]
\newtheorem{myexam}{Example}[section]

\begin{document}
\date{}
\title{ON THE CONSTRUCTION OF QUANDLES OF ORDER $3n$
\footnote{2010 mathematics subject classification primary 20N05; secondary 57M27.}
\thanks{{\bf keywords: Even quandles, examples, inner automorphism, centralizer, classification, isomorphism}}}
\author{A. O. Isere\thanks{All correspondence to be addressed to this author.}\\
Department of Mathematics,\\
Ambrose Alli University,\\
Ekpoma 310001, Nigeria.\\
abednis@yahoo.co.uk\\
isereao@aauekpoma.edu.ng \and
O. A. Elakhe\\
Department of Mathematics,\\
Ambrose Alli University,\\
Ekpoma 310001, Nigeria. \\
bssabraela@gmail.com \and
C. Ugbolo  \\
Department of Mathematics,\\
Ambrose Alli University,\\
Ekpoma 310001, Nigeria.\\
ugboloc@gmail.com }\maketitle
\begin{abstract}
 We present methods of constructing examples of quandles of order $3n, n \ge 3$. The necessary and sufficient conditions for the constructed examples to be (i)connected (ii) group (conjugate) (iii) involutory and (iv) Alexander quandles are examined and presented. Two particular examples from these methods are presented for illustration purpose and their properties are obtained, and these are used in classifying the constructed examples up to isomorphism.
\end{abstract}
\section{Introduction}
\paragraph{}
The notion of quandle was introduced independently by Joyce (\cite{qua08,qua09}) and Maveav (\cite{qua19}) in 1982 as a set with a binary operation, satisfying three axioms corresponding to Reidemeister moves of a classical knot. In knot theory, quandles play a lot of important roles, and have provided several invariants of knots. This is specially true with connected quandles (see- \cite{qua05,qua05b,qua03,qua10a,qua16,qua11}).  Strong invariants called quandle cocycles invariants were studied in \cite{qua14b} and quandles of cyclic type were studied in \cite{qua10,qua14} and were classified with cardinality up to 12. In particular, for every prime number $ p\ge 3$ there exists a quandle of cyclic type with cardinality p. This suggests that the class of quandles of cyclic type is fruitful \cite{qua10,qua14}. Moreover, the set $ C_{n}$ of Isomorphism classes of quandles of cyclic type with cardinality n were described in \cite{qua10,qua16}. Their main theorem gives a bijection from $ C_{n}$ onto $ f_{n}$ which denotes the set of cyclic permutations of order $n-1$ satisfying some conditions. This bijection is useful for studying quandles of cyclic type, since such quandles can be characterized by certain cyclic permutations. It is to be noted that all quandles of cyclic type are connected quandles.
\par
This paragraph shows the arrangement of the paper, section 2 gives a review of some relevant literature while section 3 presents the constructed quandles of order 3n, and the necessary and sufficient conditions for the constructed quandles to be connected, conjugate, involutory and Alexander quandles are presented. These are demonstrated in subsection 3.1. The problem of classifying quandles of the same order up to isomorphism, especially when they are expressed in cayley tables, is not peculiar to quandles alone \cite{qua15}. The constructed examples are classified up to isomorphism in subsection 3.2. Finally section 4 presents the concluding remarks.
\section{Preliminaries}
This section presents some definitions and results that are relevant to this work.
\begin{mydef}\label{qua06}\cite{qua05b}
A binary algebraic structure $(Q,\rhd)$ is a quandle if the following three axioms exist.
\begin{description}
\item[(1)] $ x \rhd x = x  \forall ~x\in X,$
\item[(2)] 	 There exists $z\in X $  such that $z \rhd y = x \forall ~x, y, z\in X$
\item[(3)] $ (x\rhd y) \rhd z = (x \rhd z) \rhd (y \rhd z)\forall ~x, y , z\in X$
\end{description}
\end{mydef}
\begin{mydef}\cite{qua08}\label{qua10}
A quandle $X$ equipped with two binary operations denoted as $x\rhd y$ and $x\rhd^{-1} y$ such that $x\rhd y=y^{-1}xy$
and $x\rhd^{-1} y=yxy^{-1}$ satisfying the following:
\begin{description}
\item[(1)] $ x \rhd x = x  \forall ~x\in X,$
\item[(2)] $ (x \rhd y) \rhd^{-1} y = x = (x\rhd^{-1} y)\rhd y \forall ~x, y\in X$
\item[(3)] $ (x \rhd y) * z = (x \rhd z) \rhd (y \rhd z)\forall ~x, y , z\in X$
\end{description}
is called a group (conjugate) quandle.
\end{mydef}
\begin{mydef}\label{qua09}\cite{qua06} A quandle is said to be abelian if it satisfies the identity
$(w * x)*(y * z)=(w * y)*(x * z)$.
\end{mydef}
\begin{myrem}
Abelian quandles are also medial quandles
\end{myrem}
\begin{mydef} (\cite{qua15}). Given two quandles $(X,*)$ and $(y,\rhd)$ let $f$ be a mapping from quandle
$(X,*)$  to a quandle $(Y,\rhd)$, then, $f$ is said to be a quandle homomorphism if $f(a*b)=f(a)\rhd f(b)$ for every $~x,b\in X.$
\end{mydef}
\begin{mydef} A quandle homorphism that is bijective is called quandle isomorphsim. That is to say $(X,*)$ and $(Y,\rhd)$ are isomorphic quandles if there exist an isomorphism.
\end{mydef}
The concept of isomorphism helps to distinguish two algebraic structures. Generally, the order structure has been useful
in classifying algebraic structures up to isomorphism. For some algebraic systems, the order structure is sufficient and may not be enough for others (see-\cite{qua13,qua01}).
\begin{mydef} The automorphism group of quandle $(X,*)$, denoted as $Inn(X)$ is the subgroup of $Aut(X)$ generated by all $R_{x}$ where $R_{x}(y)=y * x,$ for any $~x,y\in X$.
\end{mydef}
\begin{mydef} The inner automorphism group of a quandle $(X,*)$ denoted by $Inn(X)$, is the subgroup of $Aut(X)$ generated by all $S_{x}$, where $S_{x}(y)=y * x,$ for any $x,y\in X.$ The map $ S_{x}: X\rhd X$ that maps $u \mapsto u * x $ defines a right action of $X$ on $X$, so as to obtain a map $x \mapsto Inn(X)$.
\end{mydef}
The inner automorphism structure of quandles plays a significant role in establishing isomorphic quandles. The order of elements in this structure will to a large extent detinguish between quandles. Whenever there is a tie, then Isere (\cite{qua15}) introduced the concept of using the centralizer to establish the commutative pattern of the quandles to distinguish them.
\begin{mydef} \cite{qua15}. Centralizer of an element $'a'$ of a quandle $Q$ is the set of all members of $Q$ that commute with $'a'$.
\end{mydef}
In establishing isomorphism among quandles, you can find the centralizer of each element of the quandles one after the other until you have a different set of elements for the same element in check.
 In order to classify the constructed quandles up to isomorphism, we applied the conditions discussed by (\cite{qua13,qua15, qua01}). They said that two Quandles or loops shall be considered non-isomorphic if they contain different number of elements of the same order. Whenever, two quandles contain the same number of elements in their order structures, then we shall further consider the commutative patterns of both quandles using the centralization of certain elements in both quandles \cite{qua15}.
This is developed in subsection 3.2.

\begin{mydef} \cite{qua19c}. Let $A=Z[t^\pm1]$ be the ring of Laurent polynomials in one variable with integer coefficients. Let M be a module over A. then M is a quandle called an Alexander quandle, with quandle operation given by $a\rhd b=ta+(1-t)b.$
\end{mydef}

\begin{mydef}\label{exam1}\cite{qua05}.
Any non-empty set X with the operation
$x * y = x$ for any $~x,y\in X$ is a quandle called the trivial quandle.
\end{mydef}
\begin{mydef}\label{qua13}\cite{qua13}. A finite $Q$ is said to be connected if, for every $a,b \in Q,$ there exists $x_1 , x_2 , ..., x_n \in Q$ such that $b=x_1 \star (x_2 \star (...(...(x_n \star a))))$.

\end{mydef}
\begin{mydef}\cite{qua05} Let $n$ be a positive integer. For $a,b\in Z_{n}$  (integer modulo  n), define $a * b\equiv 2b - a$ (mod n). Then $*$ defines a quandle structure called the dihedral quandle, $R_{n}$. This set can be identified with the set of reflections of a regular n-gon with conjugation as the quandle operation.
\end{mydef}
\begin{mydef}\label{qua05}\cite{qua19d,qua08} A quandle $(X,*)$ is called involutive or involutory if $(x * y)* y = x\forall ~x,y\in X$ .
\end{mydef}
\begin{mypro}\cite{qua16, qua14} One of the following:\\
   \begin{description}
     \item[(i)] The trivial quandle $(X,s)$ is connected if and only if $ \#X = 1$
     \item[(ii)] The dihedral quandle $(X,s)$ is connected if and only if $\#X$ is odd.
     \item[(iii)] The tetrahedron quandle is connected.
   \end{description}
 \end{mypro}
  More properties of quandles of cyclic types are given (\cite{qua14,qua06,qua13}).
\begin{myexam} (\cite{qua19c}). The trivial quandle $T_{n}$ is an Alexander quandle, namely the quotient module $ T_{n} = a \rhd b = t(a) + (1-t) b = 1(a)+ (1-1) b = a.$
\end{myexam}

\begin{mypro}\label{lem1}\cite{qua05b,qua01}
  For any abelian group A, the operation $*$ defines a quandle structure on $Z_{3}\times A$ if $\mu(0) = 2,\mu(1)$ = $\mu(2) = -1,$ and $\tau(0) = 0.$ Galkin gave this definition in (\cite{qua19b}, pg. 950), for $A = Z_{p}$. The proposition generalizes his result to any abelian group A.
\begin{mylem}\cite{qua16,qua02,qua05}
   The operation $*$ as define above is idempotent and its right action is invertible  if and only if $\mu(0) = 2$ and $\tau(0) = 0$.
\end{mylem}
\end{mypro}
\begin{mypro}(\cite{qua19c}) The following are true for Alexander quandles.
\begin{enumerate}
\item Alexander Quandles are abelian .
\item If $Q$ is abelian, then $\rhd$ is left - distributive.
\end{enumerate}
\end{mypro}
\begin{myexam} (\cite{qua17})
  Let $(A,+)$ be an abelian group and $f \in Aut(A)$. Define the affine ( or Alexander) quandle over the group A as $Q_{Aff}(A,f) = (A,*),  x * y = x^{'} + y^{1-f}$.
Straight forward calculation show that $(A,*)$ is a quandle. For mediality, observe that $(x * y)* (u * v ) = (x^{'} + y^{1-f})* (u^{'} + v^{1-f} = x^{f2} + y^{(1-f)f} + u^{f(1-f)} + v^{(1-f)2}$ Is invariant under the interchange of y and u.
Alternatively, given an $R-Module$ M and an invertible element $r \in R,$ then $(M,*)$ with
$x * y = xr + y(1 - r)$
\end{myexam}
\begin{mylem}\label{qua04}
(\cite{qua19c}). If $M$ is an Alexander quandle, then for all $a, b \in M$ we have $a \rhd b + b \rhd a = a + b$
\end{mylem}
\begin{mydef}
  Let A be an abelian group. The quandle defined by $*$ in Proposition \ref{lem1} with $\mu(0) = 2, \mu(1) = \mu(2) = -1$ and $\tau(0) = 0$ is called the Galkin quandle and denoted by $G(A,\tau).$  Since $\tau $is specified by the values $\tau(1) = c_{1}$  and $\tau(2) = c_{2},$ where $c_{1},c_{2}\in A,$ we also denote it by $G(A,c_{1},c_{2})$.
\end{mydef}
\begin{mylem} For any abelian group $A$ and $c_{1} ,c_{2} \in A, G(A, c_{1},c_{2})$ and $G(A, 0,c_{2}- c_{1})$ are isomorphic.
\end{mylem}
\newpage
\begin{myexam}
One of the examples of quandles of orders 4 constructed in \cite{qua02} is presented below
\end{myexam}
\begin{table}[!hbp]
\begin{center}
$$
\begin{tabular}{|c||c|c|c||c|}
\hline
$\cdot$ & 1 & 2 & 3 & 4 \\
\hline \hline
1 & 1 & 1 & 2 & 2  \\
\hline
2 & 2 & 2 & 1 & 1 \\
\hline
3 & 4 & 4 & 3 & 3 \\
\hline
4 & 3 & 3 & 4 & 4 \\
\hline
\end{tabular}
$$
\end{center}
\caption { Quandle of Order 4}\label{qua11}
\end{table}

\section{ Main Results}
\subsection{Construction of Quandles of order 3n}
\begin{myth}\label{qua07}
 Let $(Q_{n},*)$ be any quandle of order $n$, and $n\ge 3$, and $(Z_{3},+)$ be the group of residue classes modulo 3 defined as $(Q,\rhd) = (Q_{n}\times Z_{3},\rhd)$ as follows:
 \begin{equation*}
(x,a)\rhd (y,b)=
\begin{cases}
(x*y,a+b+2)&\text{if  $a = b = 1$ or $a = 0$ and $b = 1$}\\
(x*y,a)&\text { if $a = b = 2$ or $a = 2$ and $b = 1$}\\
(x*y,a+b)&\text{if otherwise}
\end{cases}
\end{equation*} Then, $(Q,\rhd)$ is a quandle of order $3n, n\ge 3$.
\end{myth}
{\bf Proof:}\\

Let $(Q_{n},*)$ be a quandle of order n, $n\ge 3$ and  $(Z_{3},+)$ be a group of residue classes modulo 3 (With elements 0, 1 and 2).\\ We need to show that $(Q,\rhd)$ satisfies quandle axioms in Definition 2.1\\ Let $x',y',z' \in Q$ such that $x' = (x, a), y' = (y, b)$  and $z' = (z, d) ~~\forall x,y,z \in Q_{n}, a,b,d \in Z_{3}$ Axiom (1) (Idempotent).\\ First condition: For any $x' \in Q,$ then $x'\rhd x' = x'$\\ That is $(x, a)\rhd (x, a) = (x\star x, 1)=(x, a)$, where $a = b = 1$ or $a = 0$ and $b = 1$. \\
Using second condition we have, $(x, a)\rhd (x, a) = (x\star x, 2)=(x, a)$ where $a=2=b$\\ Or $(x, a)\rhd (x, a) = (x\star x, 2)=(x, a)$ where $a=2$ and $b=1$\\
Otherwise, $(x, a)\rhd (x, a) = (x\star x, 0)=(x, a)$ where $a=0=b$\\ Or $(x, a)\rhd (x, a) = (x\star x, 1)=(x, a)$ where $a=1$ and $b=0$\\
 Axiom 2. (Unique solvability). For any $y',z' \in Q$ there exists a unique element x such that $x' \rhd y' = z'$\\ That is, $(x, a) \rhd (y, b) = (z, d)$.  to show this using condition 1 where a = b = 1 or a = 0 and $b = 1$, we have $(x * y, a+b+2)=(x\star y, 1)=(z, d),$  where $z=x\star y$ and $d=1$ Or $(x\star y, 0)=(z, d)$ where $z=x\star y, d=0$\\
 Condition 2:  $(x, a) \rhd (y, b) = (x\star y, 2)=(z, d)$ where $z=x\star y, d=2$ and $a=b=2$ \\ Or $(x, a) \rhd (y, b) = (x\star y, 2)=(z, d)$ where $z=x\star y, d=2$ and $a=2$ and $b=1$ \\
 Condition 3:  $(x, a) \rhd (y, b) = (x\star y, 0)=(z, d)$ where $z=x\star y, d=0$ and $a=b=2$ \\ Or $(x, a) \rhd (y, b) = (x\star y, 1)=(z, d)$ where $z=x\star y, d=1$  \\
  Axiom 3 (Right Distributive). For any $x',y',z' \in Q$ there exist a non-associative operation such that $(x' \rhd y')\rhd z' = (x'\rhd z')\rhd (y\rhd z')$ to show that they are equal we have\\ For condition 1 $[(x, a)\rhd (y, b)]\rhd (z, d)=[(x\star y)\star z, a+b+d+1]=(x\rhd z)\rhd (y\rhd z)$.\\
  condition 2 $[(x, a)\rhd (y, b)]\rhd (z, d)=[(x\star y)\star z, a]=(x\rhd z)\rhd (y\rhd z)$ \\
  condition 3$[(x, a)\rhd (y, b)]\rhd (z, d)=[(x\star y)\star z, a+b+d]=(x\rhd z)\rhd (y\rhd z)$ \\

\begin{myth}\label{qua08} Let $(Q_{n}*)$ be a quandle of order n, $n\ge 3$ and $(Z_{3} +)$ be the group of residue classes modulo 3 defined as $(Q,\rhd) = Q_{n}\times Z_{3}$ such that
\begin{equation*}
(x,a)\rhd (y,b)=
\begin{cases}
(x*y,a+b+2)&\text{if  $a = b = 1$, $a = 0$ and $b = 1$ or $a = 2$ and $b = 1$}\\
(x*y,a+b+1)&\text { if $a = 0$ and $b = 2$ or $a = 1$ and $b = 2$}\\
(x*y,a+b)&\text{if otherwise}
\end{cases}
\end{equation*} Then, $(Q,\rhd)$ is a quandle of order $3n, n\ge 3$.
\end{myth}
{\bf Proof:}\\
The proof is similar to the proof of Theorem \ref{qua07}
\begin{myth} Let $(Q_{n},*)$ be a quandle of order n and $Z_{3}$ a residue group of integer modulo 3. Then $(\mathbb{A}=Q_{n}\times Z_{3},\rhd )$ as defined in Theorem \ref{qua07} is an involutory quandle of order 3n if and only if $(Q_n, \star)$ is an involutory quandle of order n.
\end{myth}
{\bf Proof:}\\
Let $X=(x, a)$ and $Y=(y, b)$ for all $x, y \in Q_n$ and $a, b \in Z_3$. If $(\mathbb{A}, \rhd)$ is an involutory quandle, then $(\mathbb{A}, \rhd)$ obeys definition \ref{qua05} with respect to the three conditions defined in Theorem \ref{qua07}. That is $$ (X\rhd Y)\rhd Y=X $$ implies that $$ [(x, a)\rhd (y, b)]\rhd (y, b)=(x, a) $$ implies also that $$[(x\star y)\star y, a]= (x, a) $$
$\Rightarrow (x\star y)\star y = x$. Therefore, $(Q_n, \star)$ is an involutory quandle.\\
Conversely, suppose $(Q_n, \star)$ is an involutory quandle, then we need to show that $(\mathbb{A}, \rhd)$ as defined in Theorem \ref{qua07} is an involutory quandle. That is, $(\mathbb{A}, \rhd)$ obeys definition 2.2\\
Condition 1: consider
$$
(X\rhd Y)\rhd Y=[(x, a)\rhd (y, b)]\rhd (y, b)=(x\star y, a+b+2)\rhd (y, b)=[(x\star y)\star y, a+2b+1]=[(x\star y)\star y, a]
$$ where $a=b=1$ or $a=0$ and $b=1$.\\
Since $(Q_n, \star)$ is an involutory quandle, then $$ [(x\star y)\star y, a]=(x, a)=X $$
Condition 2: $$ (X\rhd Y)\rhd Y =[(x, a)\rhd (y, b)]\rhd (y, b)=(x\star y, a)\rhd (y, b)=[(x\star y)\star y, a]$$ where $a=b=2$ or $a=2$ and $b=1$.\\
Since $(Q_n, \star)$ is an involutory quandle, then $$ [(x\star y)\star y, a]=(x, a)=X $$
Condition 3:
$$ (X\rhd Y)\rhd Y =[(x, a)\rhd (y, b)]\rhd (y, b)=(x\star y, a+b)\rhd (y, b)=[(x\star y)\star y, a+2b]=[(x\star y)\star y, a]$$ where $a=b=0$ or $a=1$ and $b=0$.\\
Since $(Q_n, \star)$ is an involutory quandle, then $$ [(x\star y)\star y, a]=(x, a)=X $$
Therefore, $$ (X\rhd Y)\rhd Y=X .$$ Thus, $(\mathbb{A}, \rhd)$ is an involutory quandle.
\begin{myrem}
That the quandle $(\mathbb{B}=Q_n \times Z_3 , \rhd)$ as defined in Theorem \ref{qua08} is an involutory quandle is proved in a similar manner.
\end{myrem}

\begin{myth} Let $(Q_{n},\star)$ be a quandle of order n and $Z_{3}$ a residue group of integer modulo 3. Then $(\mathbb{A}=Q_{n}\times Z_{3},\rhd )$ as defined in Theorem \ref{qua07} is a conjugate (group) quandle of order 3n if and only if $(Q_n, \star)$ is a conjugate (group) quandle of order n.
\end{myth}
{\bf Proof:}\\
Let $X=(x, a)$ and $Y=(y, b)$ for all $x, y \in Q_n$ and $a, b \in Z_3$, and if $(\mathbb{A}, \rhd)$ is a group quandle, then $(\mathbb{A}, \rhd)$ obeys definition \ref{qua10} with respect to the three conditions defined in Theorem \ref{qua07}. That is $$ (X\rhd Y)=Y^{-1}XY $$ and $$ (X\rhd^{-1} Y)=YXY^{-1} $$ First,  $$ X\rhd Y=(y, b)^{-1}(x, a)(y, b)= (y^{-1}, -b)(x, a)(y, b) $$ implies that $$(x, a)\rhd (y, b)= (y^{-1}xy, a). $$ That is $$(x\star y, a)= (y^{-1}xy, a). $$
Secondly, $$ X\rhd^{-1} Y=(y, b)(x, a)(y, b)^{-1}= (y, b)(x, a)(y^{-1}, -b) $$ implies that $$ (x, a)\rhd^{-1} (y, b)=(x\star y, a)= (yxy^{-1}, a). $$ Therefore,  $$x\star y= y^{-1}xy ~~\&~~ x\star y=yxy^{-1}$$ Therefore, $(Q_{n},\star)$ is a conjugate (group) quandle.\\
Conversely, suppose $(Q_{n},\star)$ is a conjugate (group) quandle, then we need to show that $(\mathbb{A}, \rhd )$ as defined in Theorem \ref{qua07} is a conjugate (group) quandle. That is that $(\mathbb{A}, \rhd)$ obeys definition 2.2.\\

condition 1:

$$
(X\rhd Y)=[(x, a)\rhd (y, b)]=(x\star y, a+b+2)=(x\star y, a)=(y^{-1}xy, a)=Y^{-1}XY
$$ where $a=b=1$ or $a=0$ and $b=1$.\\ and $$
(X\rhd^{-1} Y)=[(x, a)\rhd (y, b)]=(x\star y, a)=(yxy^{-1}, a)=YXY^{-1}
$$ where $a=b=1$ or $a=0$ and $b=1$.\\
Since $(Q_n, \star)$ is a conjugate (group) quandle, where $$ x\star y =y^{-1}xy ~~\&~~ x\star y = yxy^{-1} $$ respectively.\\
Therefore, $X\rhd Y = Y^{-1}XY$ and $X\rhd^{-1} Y=YXY^{-1}$\\
Condition 2: $$ X\rhd Y =(x\star y, a)=(y^{-1}xy, a)=Y^{-1}XY~~\&~~ X\rhd^{-1} Y =(x\star y, a)=(y^{-1}xy, a)=YXY^{-1}$$ where $a=b=2$ or $a=2$ and $b=1$.\\
Condition 3:
$$ X\rhd Y =(x\star y, a+b)=(y^{-1}xy, a)=Y^{-1}XY~~\&~~ X\rhd^{-1} Y =(x\star y, a+b)=(yxy^{-1}, a)==YXY^{-1}.$$
where $a=b=0$ or $a=1$ and $b=0$.\\

Therefore, $$X\rhd Y=Y^{-1}XY $$ and $$ X\rhd^{-1} Y=YXY^{-1}.$$  The proof is complete.

\begin{mylem}\label{qua02} Let $(Q_{n},\star)$ be a quandle of order n and $Z_{3}$ a residue group of integer modulo 3. Then $(Q =Q_{n}\times Z_{3},\rhd )$ as defined in Theorem \ref{qua07} is a left distributive quandle if and only if $(Q_n, \star)$ is left distributive.
\end{mylem}
{\bf Proof:}\\
Let $X=(x, a), Y=(y, b) $ and $Z=(z, c)$ for all $x, y, z \in Q_n$ and $a, b, c \in Z_3$, and if $(Q, \rhd)$ is a left distributive quandle, then $(Q, \rhd)$ obeys $X \rhd (Y \rhd Z)=(X \rhd Y) \rhd (X \rhd Z)$ with respect to the three conditions defined in Theorem \ref{qua07}. That is $$ X \rhd (Y \rhd Z)=[x\star (y\star z), a] $$ and $$ (X \rhd Y) \rhd (X \rhd Z)= [(x\star y)\star (x\star z), a].$$ Thus, $$ x\star (y\star z)= (x\star y)\star (x\star z) .$$ Therefore, $(Q_n, \star)$ is a left distributive quandle.\\
Conversely, given that $(Q_n, \star)$ is left distributive, we want to show that $(Q, \rhd)$ is also left distributive subject to the three sets of conditions.\\
Condition 1: Consider $$ X \rhd (Y \rhd Z)= (x, a)\rhd [(y, b)\rhd (z, c)]=(x, a)\rhd (y\star z, b+c+2)=[x\star (y\star z), a+b+c+1]=[x\star (y\star z), a]$$ where $a=b=c=1$ or $a=0$ and $b=c=1$. Since $(Q_n, \star)$ is left distributive, then $$[x\star (y\star z), a]= [(x\star y)\star (x\star z), a]=(X \rhd Y) \rhd (X \rhd Z).$$
Condition 2:  $$ X \rhd (Y \rhd Z)= (x, a)\rhd [(y, b)\rhd (z, c)]=(x, a)\rhd (y\star z, b)=[x\star (y\star z), a]$$ where $a=b=c=2$ or $a=2$ and $b=c=1$. Implies that $$[x\star (y\star z), a]= [(x\star y)\star (x\star z), a]=(X \rhd Y) \rhd (X \rhd Z).$$ Since $(Q_n, \star)$ is left distributive.\\
Condition 3:  $$ X \rhd (Y \rhd Z)= (x, a)\rhd [(y, b)\rhd (z, c)]=(x, a)\rhd (y\star z, b+c)=[x\star (y\star z), a+b+c]=[x\star (y\star z), a]$$ where $a=b=c=0$ or $a=1$ and $b=c=0$. Implies that $$[x\star (y\star z), a]= [(x\star y)\star (x\star z), a]=(X \rhd Y) \rhd (X \rhd Z).$$ Since $(Q_n, \star)$ is left distributive. Therefore, $(Q, \rhd)$ is a left distributive quandle.

\begin{mylem}\label{qua03} Let $(Q_{n},\star)$ be a quandle of order n and $Z_{3}$ a residue group of integer modulo 3. Then $(Q =Q_{n}\times Z_{3},\rhd )$ as defined in Theorem \ref{qua07} is abelian if and only if $(Q_n, \star)$ is abelian.
\end{mylem}
{\bf Proof:}\\
Let $A=(x, a), B=(y, b), C=(z, c) $ and $D=(t, d)$ for all $x, y, z, t \in Q_n$ and $a, b, c, d \in Z_3$, and if $(Q, \rhd)$ is abelian, then $(Q, \rhd)$ obeys definition \ref{qua09} That is $(A \rhd B)\rhd (C \rhd D)=(A \rhd C) \rhd (B\rhd D)$ with respect to the three conditions defined in Theorem \ref{qua07}. Consider: $$ (A \rhd B)\rhd (C \rhd D)=[(x, a)\rhd (y, b)]\rhd[(z, c)\rhd (t, d)] = [(x\star y )\star (z\star t), a] $$ and $$ (A \rhd C) \rhd (B\rhd D)= [(x\star z)\star (y\star t), a].$$ Thus,  $$ (x\star y)\star (z\star t)= (x\star z)\star (y\star t) \forall x, y, z, t \in Q_n$$ Therefore, $(Q_n, \star)$ is abelian.\\
Conversely, given that $(Q_n, \star)$ is abelian, we want to show that $(Q, \rhd)$ is also abelian subject to the three sets of conditions defined in Theorem \ref{qua07}.\\
Condition 1: $$ (A \rhd B)\rhd (C \rhd D)= [(x, a)\rhd (y, b)]\rhd[(z, c)\rhd (t, d)] =(x\star y, a+b+2)\rhd (z\star t, c+d+2)=[(x\star y)\star (z\star t), a]$$ where $a=b=c=d=1$ or $a=0$ and $b=c=d=1$. Since $(Q_n, \star)$ is abelian, then $$[(x\star y)\star (z\star t), a]= [(x\star z)\star (y\star t), a]=(A \rhd C) \rhd (B\rhd D).$$
Condition 2:  $$ (A \rhd B)\rhd (C \rhd D)= [(x, a)\rhd (y, b)]\rhd[(z, c)\rhd (t, d)] =(x\star y, a)\rhd (z\star t, c)=[(x\star y)\star (z\star t), a].$$  where $a=b=c=d=2$ or $a=2$ and $b=c=d=1$. Since $(Q_n, \star)$ is abelian, then $$[(x\star y)\star (z\star t), a]= [(x\star z)\star (y\star t), a]=(A \rhd C) \rhd (B\rhd D).$$
Condition 3: Otherwise,  $$ (A \rhd B)\rhd (C \rhd D)= [(x\star y)\star (z\star t), a].$$  where $a=b=c=d=1$ or $a=2$ and $b=c=d=0$. Since $(Q_n, \star)$ is abelian, then $$[(x\star y)\star (z\star t), a]= [(x\star z)\star (y\star t), a]=(A \rhd C) \rhd (B\rhd D).$$ Therefore, $(Q, \rhd)$ is a abelian.

\begin{myth} Let $(Q_{n},\star)$ be a quandle of order n and $Z_{3}$ a residue group of integer modulo 3. Then $(Q=Q_{n}\times Z_{3},\rhd )$ as defined in Theorem \ref{qua07} is an Alexander quandle of order 3n if and only if $(Q_n, \star)$ is an Alexander quandle of order n.
\end{myth}
{\bf Proof:}\\
Let $X=(x, a)$ and $Y=(y, b)$ for all $x, y \in Q_n$ and $a, b \in Z_3$, and if $(Q, \rhd)$ is an Alexander quandle, then $(Q, \rhd)$ obeys Lemma \ref{qua04} with respect to the three conditions defined in Theorem \ref{qua07}. That is $$ X\rhd Y + Y\rhd X = X + Y .$$ This implies that
$$ X\rhd Y + Y\rhd X = [(x, a)\rhd (y, b)]+[(y, b)\rhd (x, a)]=(x\star y, a)+(y\star x, b)=(x\star y+y\star x, a+b)$$ and
$$ X+Y=(x, a)+(y, b)=(x+y, a+b).$$ Implies that $$ x\star y + y\star x=x+y. $$ Therefore, $(Q_n, \star)$ obeys Lemma \ref{qua04}
and thus, it is an Alexander quandle.\\
Conversely, suppose that $(Q_n, \star)$ is an Alexander quandle then, we need to show that $(Q, \rhd)$ as defined in Theorem \ref{qua07} is also an Alexander quandle. That $(Q, \rhd)$ obeys Lemma \ref{qua07}.\\
Condition 1: consider $$ X\rhd Y + Y\rhd X = (x, a)\rhd (y, b)+(y, b)\rhd (x, a)=(x\star y, a)+(y\star x, b)=(x\star y+y\star x, a+b)$$ where $a=b=1$ or $a=0$ ~~\&~~ $b=1$ and $b=a=1$ or $b=0$ ~~\&~~$a=1$ respectively. Implies that $$ (x\star y+y\star x, a+b)=(x+y, a+b)=X + Y $$ Since, $(Q_n, \star)$ obeys Lemma \ref{qua07}.\\
Condition 2: consider $$ X\rhd Y + Y\rhd X = (x, a)\rhd (y, b)+(y, b)\rhd (x, a)=(x\star y, a)+(y\star x, b)=(x\star y+y\star x, a+b)$$ where $a=b=2$ or $a=2$ ~~\&~~ $b=1$ and $b=a=2$ or $b=2$ ~~\&~~$a=1$ respectively. Implies that $$ (x\star y+y\star x, a+b)=(x+y, a+b)=X + Y $$ Since, $(Q_n, \star)$ obeys Lemma \ref{qua07}.\\
Condition 3: consider $$ X\rhd Y + Y\rhd X = (x, a)\rhd (y, b)+(y, b)\rhd (x, a)=(x\star y, a+b)+(y\star x, b+a)=(x\star y+y\star x, a+b)$$ where $a=b=0$ or $a=1$ ~~\&~~ $b=0$ and $b=a=0$ or $b=1$ ~~\&~~$a=0$. Implies that $$ (x\star y+y\star x, a+b)=(x+y, a+b)=X+Y $$ Since, $(Q_n, \star)$ obeys Lemma \ref{qua07}.\\
Therefore, $$ X\rhd Y+Y\rhd X=X+Y.$$ Thus, $(Q, \rhd)$ is an Alexander quandle

\begin{myth} Let $(Q_{n},\star)$ be a quandle of order n and $Z_{3}$ a residue group of integer modulo 3. Then $(Q=Q_{n}\times Z_{3},\rhd )$ as defined in Theorem \ref{qua07} is a connected quandle of order 3n if and only if $(Q_n, \star)$ is also connected.
\end{myth}
{\bf Proof:}\\
Let $A=(x, a), B=(y, b), X=(x_{1}, a_{1})$ and $Y=(y_{1}, b_{1})$ for all $x, y,x_{1}, y_{1} \in Q_n$ and $a, b, a_{1}, b_{1}\in Z_3$, and if $(Q, \rhd)$ is a connected quandle, then $(Q, \rhd)$ obeys Definition \ref{qua13} with respect to the three conditions defined in Theorem \ref{qua07}. That is $$ X\rhd (Y\rhd A)=B$$ implies that $$ (x_{1}, a_{1})\rhd [(y_{1}, b_{1})\rhd (x, a)]=[x_{1}\star(y_{1}\star x), a_{1}]=(y, b).$$ Then, $$ x_{1}\star(y_{1}\star x)=y.$$ Therefore, $(Q_n, \star)$ is a connected quandle whenever $a_{1}=b$.\\
Conversely, suppose $(Q_n, \star)$ is a connected quandle, then we need to show $(Q, \rhd)$ is also connected.\\
Condition 1:$$ X\rhd (Y\rhd A)= (x_{1}, a_{1})\rhd [(y_{1}, b_{1})\rhd (x, a)]=[x_{1}\star(y_{1}\star x), a_{1}+b_{1}+2]=[x_{1}\star(y_{1}\star x), a_{1}]= (y, a_{1})$$ where $a_{1}=b_{1}=1$ or $a_{1}=0$ and $b_{1}=1$.\\
Condition 2:$$ X\rhd (Y\rhd A)= (x_{1}, a_{1})\rhd [(y_{1}, b_{1})\rhd (x, a)]=[x_{1}\star(y_{1}\star x), a_{1}]=(y, a_{1}$$ where $a_{1}=b_{1}=2$ or $a_{1}=2$ and $b_{1}=1$.\\
Condition 3:$$ X\rhd (Y\rhd A)= (x_{1}, a_{1})\rhd [(y_{1}, b_{1})\rhd (x, a)]=[x_{1}\star(y_{1}\star x), a_{1}+a+b_{1}]=[x_{1}\star(y_{1}\star x), a_{1}]=(y, a_{1}$$ where $a_{1}=a=b_{1}=0$ or $a_{1}=1$ and $a=b_{1}=0$.\\
Therefore $$[x_{1}\star(y_{1}\star x), a_{1}]=(y, b) $$ whenever $a_{1}=b$.

\newpage
\subsection{Classification of the Constructed Examples up to Isomorphism}

\begin{myexam}\label{qua01a}
A quandle of order 12 constructed from a quandle of order 4 in Table \ref{qua11} using Theorem \ref{qua07}.
\end{myexam}
\begin{table}[!hbp]
\begin{center}
\begin{tabular}{|c|c|c|c|c|c|c|c|c|c|c|c|c|}
\hline
$\rhd$ & 1 & 2 & 3 & 4 & 5 & 6 & 7 & 8 & 9 & 10 & 11 & 12 \\
\hline \hline
1 & 1 & 1 & 1 &	1 &	1 &	1 &	1 &	1 &	1 &	1 &	1 &	1 \\
\hline
2 & 2 & 2 & 2 & 2 &	2 &	2 &	2 &	2 &	2 &	2 &	2 &	2 \\
\hline
3 &	3 &	3 &	3 &	3 &	3 &	3 &	3 &	3 &	3 &	3 &	3 &	3 \\
\hline
4 &	4 &	4 &	4 &	4 &	4 &	4 &	10 & 10 & 10 &	7 &	7 &	7 \\
\hline
5 &	5 &	5 &	5 &	5 &	5 &	5 &	11 & 11 & 11 &	8 &	8 &	8 \\
\hline
6 &	6 &	6 &	6 &	6 &	6 &	6 &	12 & 12 & 12 &	9 &	9 &	9 \\
\hline
7 &	7 & 7 & 7 & 10 & 10 & 10 & 7 & 7 & 7 & 4 & 4 & 4 \\
\hline
8 &	8 & 8 & 8 & 11 & 11 & 11 & 8 & 8 & 8 & 5 & 5 & 5 \\
\hline
9 & 9 & 9 & 9 & 12 & 12 & 12 & 9 & 9 & 9 & 6 & 6 & 6 \\
\hline
10 & 10 & 10 & 10 & 7 & 7 & 7 & 4 & 4 & 4 & 10 & 10 & 10 \\
\hline
11 & 11 & 11 & 11 & 8 & 8 & 8 & 5 & 5 & 5 & 11 & 11 & 11 \\
\hline
12 & 12 & 12 & 12 & 9 & 9 & 9 & 6 & 6 & 6 & 12 & 12 & 12 \\
\hline
\end{tabular}
\end{center}
\caption{$Q_{1}$ - A Quandle of Order 12 By Theorem \ref{qua07}}\label{qua12a}
\end{table}
 {\bf Inner Automorphism Structure of $Q_{1}$}\\
R(1) = (1)\\
R(2) = (1)\\
R(3) = (1)\\
R(4) = (7,10), (8,11), (9,12)\\
R(5) = (7,10), (8,11), (9,12)\\
R(6) = (7,10), (8,11), (9,12)\\
R(7) = (4,10), (5,11), (6,12)\\
R(8) = (4,10), (5,11), (6,12)\\
R(9) = (4,10), (5,11), (6,12)\\
R(10)= (4,7),  (5,8),  (6,9)\\
R(11)= (4,7),  (5,8),  (6,9)\\
R(12)= (4,7),  (5,8),  (6,9)\\
The inner automorphism structure above shows that this quandle ($Q_{1}$) of order 12 in Table \ref{qua12a} has 9 elements of order 2.
\newpage
 \begin{myexam}\label{qua01b}
Another  quandle of order 12 constructed from the same quandle of order 4 in Table \ref{qua11} using Theorem \ref{qua07}.
\end{myexam}
\begin{table}[!hbp]
\begin{center}
\begin{tabular}{|c|c|c|c|c|c|c|c|c|c|c|c|c|}
\hline
$\rhd$ & 1 & 2 & 3 & 4 & 5 & 6 & 7 & 8 & 9 & 10 & 11 & 12 \\
\hline \hline
1 & 1 & 1 & 2 & 1 & 1 & 2 & 1 & 1 & 2 & 1 & 1 & 2 \\
\hline
2 & 2 & 2 & 1 & 2 & 2 & 1 & 2 & 2 & 1 & 2 & 2 & 1 \\
\hline
3 & 3 & 3 & 3 & 3 & 3 & 3 & 3 & 3 & 3 & 3 & 3 & 3 \\
\hline
4 & 4 & 4 & 5 & 4 & 4 & 5 & 10 & 10 & 11 & 7 & 7 & 8 \\
\hline
5 & 5 & 5 & 4 & 5 & 5 & 4 & 11 & 11 & 10 & 8 & 8 & 7 \\
\hline
6 & 6 & 6 & 6 & 6 & 6 & 6 & 12 & 12 & 12 & 12 & 12 & 12 \\
\hline
7 & 7 & 7 & 8 & 10 & 10 & 11 & 7 & 7 & 8 & 4 & 4 & 5 \\
\hline
8 & 8 & 8 & 7 & 11 & 11 & 10 & 8 & 8 & 7 & 5 & 5 & 4 \\
\hline
9 & 9 & 9 & 9 & 12 & 12 & 12 & 9 & 9 & 9 & 6 & 6 & 6 \\
\hline
10 & 10 & 10 & 11 & 7 & 7 & 8 & 4 & 4 & 5 & 10 & 10 & 11 \\
\hline
11 & 11 & 11 & 10 & 8 & 8 & 7 & 5 & 5 & 4 & 11 & 11 & 10 \\
\hline
12 & 12 & 12 & 12 & 9 & 9 & 9 & 6 & 6 & 6 & 12 & 12 & 12 \\
\hline
\end{tabular}
\end{center}
\caption{$Q_{2}$ - A Quandle of Order 12 by Theorem \ref{qua08}}\label{qua12}
\end{table}
{\bf Inner Automorphism Structure of $Q_{2}$}\\
R(1) = (1)\\
R(2) = (1)\\
R(3) = (1,2), (4,5), (7,8), (10,11)\\
R(4) = (7,10), (8,11), (9,12)\\
R(5) = (7,10), (8,11), (9,12)\\
R(6) = (1,2), (4,5), (7,11), (8,10), (9,12)\\
R(7) = (4,10), (5,11), (6,12)\\
R(8) = (4,10), (5,11), (6,12)\\
R(9) = (1,2), (4,11), (5,10), (6,12)\\
R(10) = (4,7), (5,8), (6,9)\\
R(11) = (4,7), (5,8), (6,9)\\
R(12) = (1,2), (4,8), (5,7), (6,9)\\
The inner automorphism structure shows that this quandle ($Q_{2}$) of order 12 in Table \ref{qua12} has 10 elements of order 2.

\begin{myth} The constructed examples of quandles of order 3n in Theorem \ref{qua07} and Theorem \ref{qua08} are non-isomorphic.
\end{myth}
{\bf Proof:}\\
The proof follows from Example \ref{qua01a} and Example \ref{qua01b}.

\section{Concluding Remarks}
The paper discussed two methods of constructing quandles of order 3n. The results obtained for Theorem \ref{qua07} can also be proved for Theorem \ref{qua08} in a similar manner. The classification of the constructed examples up to isomorphism was demonstrated in subsection 3.2 using two illustrative examples wherein their inner automorphism structures were sufficient to distinguish these quandles up to isomorphism. It is to be noted however that not in all cases that the inner automorphism structures of quandles will be sufficient for classification of quandles up to isomorphism (see-\cite{qua15}).

\end{document}